\documentstyle[twoside]{article}

\newread\epsffilein    % file to \read
\newif\ifepsffileok    % continue looking for the bounding box?
\newif\ifepsfbbfound   % success?
\newif\ifepsfverbose   % report what you're making?
\newif\ifepsfdraft     % use draft mode?
\newdimen\epsfxsize    % horizontal size after scaling
\newdimen\epsfysize    % vertical size after scaling
\newdimen\epsftsize    % horizontal size before scaling
\newdimen\epsfrsize    % vertical size before scaling
\newdimen\epsftmp      % register for arithmetic manipulation
\newdimen\pspoints     % conversion factor
\def\epsfbox#1{\global\def\epsfllx{72}\global\def\epsflly{72}%
   \global\def\epsfurx{540}\global\def\epsfury{720}%
   \def\lbracket{[}\def\testit{#1}\ifx\testit\lbracket
   \let\next=\epsfgetlitbb\else\let\next=\epsfnormal\fi\next{#1}}%
\def\epsfgetlitbb#1#2 #3 #4 #5]#6{\epsfgrab #2 #3 #4 #5 .\\%
   \epsfsetgraph{#6}}%
\def\epsfnormal#1{\epsfgetbb{#1}\epsfsetgraph{#1}}%
\def\epsfgetbb#1{%
\openin\epsffilein=#1
\ifeof\epsffilein\errmessage{I couldn't open #1, will ignore it}\else
   {\epsffileoktrue \chardef\other=12
    \def\do##1{\catcode`##1=\other}\dospecials \catcode`\ =10
    \loop
       \read\epsffilein to \epsffileline
       \ifeof\epsffilein\epsffileokfalse\else
          \expandafter\epsfaux\epsffileline:. \\%
       \fi
   \ifepsffileok\repeat
   \ifepsfbbfound\else
    \ifepsfverbose\message{No bounding box comment in #1; using defaults}\fi\fi
   }\closein\epsffilein\fi}%
\def\epsfclipoff{\def\epsfclipstring{\ifepsfdraft\space clip\fi}}%
\epsfclipoff
\def\epsfsetgraph#1{%
   \epsfrsize=\epsfury\pspoints
   \advance\epsfrsize by-\epsflly\pspoints
   \epsftsize=\epsfurx\pspoints
   \advance\epsftsize by-\epsfllx\pspoints
   \epsfxsize\epsfsize\epsftsize\epsfrsize
   \ifnum\epsfxsize=0 \ifnum\epsfysize=0
      \epsfxsize=\epsftsize \epsfysize=\epsfrsize
      \epsfrsize=0pt
     \else\epsftmp=\epsftsize \divide\epsftmp\epsfrsize
       \epsfxsize=\epsfysize \multiply\epsfxsize\epsftmp
       \multiply\epsftmp\epsfrsize \advance\epsftsize-\epsftmp
       \epsftmp=\epsfysize
       \loop \advance\epsftsize\epsftsize \divide\epsftmp 2
       \ifnum\epsftmp>0
          \ifnum\epsftsize<\epsfrsize\else
             \advance\epsftsize-\epsfrsize \advance\epsfxsize\epsftmp \fi
       \repeat
       \epsfrsize=0pt
     \fi
   \else \ifnum\epsfysize=0
     \epsftmp=\epsfrsize \divide\epsftmp\epsftsize
     \epsfysize=\epsfxsize \multiply\epsfysize\epsftmp   
     \multiply\epsftmp\epsftsize \advance\epsfrsize-\epsftmp
     \epsftmp=\epsfxsize
     \loop \advance\epsfrsize\epsfrsize \divide\epsftmp 2
     \ifnum\epsftmp>0
        \ifnum\epsfrsize<\epsftsize\else
           \advance\epsfrsize-\epsftsize \advance\epsfysize\epsftmp \fi
     \repeat
     \epsfrsize=0pt
    \else
     \epsfrsize=\epsfysize
    \fi
   \fi
   \ifepsfverbose\message{#1: width=\the\epsfxsize, height=\the\epsfysize}\fi
   \epsftmp=10\epsfxsize \divide\epsftmp\pspoints
   \vbox to\epsfysize{\vfil\hbox to\epsfxsize{%
      \ifnum\epsfrsize=0\relax
        \includegraphics{\ifepsfdraft}%
      \else
        \epsfrsize=10\epsfysize \divide\epsfrsize\pspoints
        \includegraphics{\ifepsfdraft}%
      \fi
      \hfil}}%
\global\epsfxsize=0pt\global\epsfysize=0pt}%
{\catcode`\%=12 \global\let\epsfpercent=%\global\def\epsfbblit{%BoundingBox}}%
\long\def\epsfaux#1#2:#3\\{\ifx#1\epsfpercent
   \def\testit{#2}\ifx\testit\epsfbblit
      \epsfgrab #3 . . . \\%
      \epsffileokfalse
      \global\epsfbbfoundtrue
   \fi\else\ifx#1\par\else\epsffileokfalse\fi\fi}%
\def\epsfempty{}%
\def\epsfgrab #1 #2 #3 #4 #5\\{%
\global\def\epsfllx{#1}\ifx\epsfllx\epsfempty
      \epsfgrab #2 #3 #4 #5 .\\\else
   \global\def\epsflly{#2}%
   \global\def\epsfurx{#3}\global\def\epsfury{#4}\fi}%
\def\epsfsize#1#2{\epsfxsize}
\catcode`\@=11
\long\def\@makefntext#1{
\protect\noindent \hbox to 3.2pt {\hskip-.9pt
$^{{\eightrm\@thefnmark}}$\hfil}#1\hfill}		%CAN BE USED

\def\@makefnmark{\hbox to 0pt{$^{\@thefnmark}$\hss}}	%ORIGINAL
	
\def\ps@myheadings{\let\@mkboth\@gobbletwo
\def\@oddhead{\hbox{}
\rightmark\hfil\eightrm\thepage}
\def\@oddfoot{}\def\@evenhead{\eightrm\thepage\hfil
\leftmark\hbox{}}\def\@evenfoot{}
\def\sectionmark##1{}\def\subsectionmark##1{}}

\oddsidemargin=\evensidemargin
\addtolength{\oddsidemargin}{-30pt}
\addtolength{\evensidemargin}{-30pt}

\newcounter{sectionc}\newcounter{subsectionc}\newcounter{subsubsectionc}
\renewcommand{\section}[1] {\vspace{12pt}\addtocounter{sectionc}{1}
\setcounter{subsectionc}{0}\setcounter{subsubsectionc}{0}\noindent
	{\tenbf\thesectionc. #1}\par\vspace{5pt}}
\renewcommand{\subsection}[1] {\vspace{12pt}\addtocounter{subsectionc}{1}
	\setcounter{subsubsectionc}{0}\noindent
	{\bf\thesectionc.\thesubsectionc. #1}\par\vspace{5pt}}
\renewcommand{\subsubsection}[1] {\vspace{12pt}\addtocounter{subsubsectionc}{1}
        \noindent{\bf\thesectionc.\thesubsectionc.\thesubsubsectionc.
	#1}\par\vspace{5pt}}

\topsep=0in\parsep=0in\itemsep=0in
\parindent=15pt

\newcommand{\textlineskip}{\baselineskip=13pt}
\newcommand{\smalllineskip}{\baselineskip=10pt}

\def\eightcirc{
\begin{picture}(0,0)
\put(4.4,1.8){\circle{6.5}}
\end{picture}}
\def\eightcopyright{\eightcirc\kern2.7pt\hbox{\eightrm c}}

\def\abstracts#1#2#3{{
	\centering{\begin{minipage}{4.5in}\baselineskip=10pt\footnotesize
        \parindent=0pt {\centering{ABSTRACT\par\vspace{2pt}}}
	\parindent=15pt #1\par  %% better 0pt
	\parindent=15pt #2\par
	\parindent=15pt #3
	\end{minipage}}\par}}

\def\keywords#1{{
	\centering{\begin{minipage}{4.5in}\baselineskip=10pt\footnotesize
	{\footnotesize\it Keywords}\/: #1
	 \end{minipage}}\par}}

\newcounter{itemlistc}
\newcounter{romanlistc}
\newcounter{alphlistc}
\newcounter{arabiclistc}

\newcommand{\fcaption}[1]{
        \refstepcounter{figure}
        \setbox\@tempboxa = \hbox{\footnotesize Fig.~\thefigure. #1}
        \ifdim \wd\@tempboxa > 5in
           {\begin{center}
        \parbox{5in}{\footnotesize\smalllineskip Fig.~\thefigure. #1}
            \end{center}}
        \else
             {\begin{center}
             {\footnotesize Fig.~\thefigure. #1}
              \end{center}}
        \fi}

\newcommand{\tcaption}[1]{
        \refstepcounter{table}
        \setbox\@tempboxa = \hbox{\footnotesize Table~\thetable. #1}
        \ifdim \wd\@tempboxa > 5in
           {\begin{center}
        \parbox{5in}{\footnotesize\smalllineskip Table~\thetable. #1}
            \end{center}}
        \else
             {\begin{center}
             {\footnotesize Table~\thetable. #1}
              \end{center}}
        \fi}

\def\fnt#1#2{\footnotetext{\kern-.3em
   {$^{\mbox{\scriptsize #1}}$}{#2}}}

\def\runninghead#1#2{\pagestyle{myheadings}
\markboth{{\protect\footnotesize\it{\quad #1}}\hfill}
{\hfill{\protect\footnotesize\it{#2\quad}}}}
\headsep=15pt

\font\tenit=cmti10
\font\tenbf=cmbx10

\font\ninerm=cmr9
\font\nineit=cmti9
\font\ninebf=cmbx9
\font\eightrm=cmr8

\textwidth=5truein
\textheight=7.8truein

\def\qed{\hbox{${\vcenter{\vbox{			%HOLLOW SQUARE
   \hrule height 0.4pt\hbox{\vrule width 0.4pt height 6pt
   \kern5pt\vrule width 0.4pt}\hrule height 0.4pt}}}$}}

\textheight=7.8truein
\textwidth=5truein

\addtolength{\oddsidemargin}{-0.3cm}   %%% changed here
\addtolength{\evensidemargin}{-0.3cm}  %%% and here
\setcounter{footnote}{0}

\newcommand\firstline{\noindent}
\newcommand\qedp{. \qed}
\newcommand\skipline{\par\vspace*{12pt}\noindent}
\newcommand\skipsline{\par\vspace*{6pt}\noindent}

\newcommand\fighead{\vspace*{20pt}}
\newcommand\figtail[1]{\vspace*{6pt}\fcaption{#1}\vspace*{13pt}}

\newcommand\showfbox{}

\newcounter{definitionc}
\setcounter{definitionc}{0}
\newenvironment{definition}{\refstepcounter{definitionc}\bf Definition \thedefinitionc. \rm}{}

\newcounter{theoremc}
\setcounter{theoremc}{0}
\newenvironment{theorem}{\refstepcounter{theoremc}\bf Theorem \thetheoremc. \it}{\rm}

\newcounter{corollaryc}
\setcounter{corollaryc}{0}
\newenvironment{corollary}{\refstepcounter{corollaryc}\bf Corollary \thecorollaryc. \it}{\rm}

\newcounter{lemmac}
\setcounter{lemmac}{0}
\newenvironment{lemma}{\refstepcounter{lemmac}\bf Lemma \thesectionc.\thelemmac. \it}{\rm}

\newcounter{mylistc}
\newenvironment{mylist}[1]
	{\setcounter{mylistc}{0}
	 \begin{list}{$($\alph{mylistc}$)$}
	{\usecounter{mylistc}
	 \setlength{\listparindent}{\parindent}
	 \setlength{\partopsep}{0pt}
	 \setlength{\parsep}{0pt}
	 \setlength{\topsep}{1pt}
	 \setlength{\itemsep}{1pt}#1}}{\end{list}}

\newcommand\proof[1]{{\noindent}{\bf Proof#1. }}
\newcommand\remark{{\noindent}{\bf Remark. }}

\newcommand\Bn{\ifmmode {\cal B}_n\else ${\cal B}_n$\fi}
\newcommand\Pn{\ifmmode {\cal P}_{\!n}\else ${\cal P}_{\!n}$\fi}
\newcommand\Wn{\ifmmode {\cal W}_n\else ${\cal W}_n$\fi}
\newcommand\Wnk{\ifmmode {\cal W}_{n}^{k}\else ${\cal W}_{n}^{k}$\fi}
\newcommand\Wne{\ifmmode {\cal W}_{n}^{1}\else ${\cal W}_{n}^{1}$\fi}
\newcommand\Kn{\ifmmode {\cal K}_n\else ${\cal K}_n$\fi}
\newcommand\squad{\em}
\newcommand\pink{\ifmmode \pi_{n_{1}\ldots n_{k}}\else $\pi_{n_{1} n_{2} \ldots n_{k}}$\fi}
\newcommand\mathpink{\pi_{n_{1}\ldots n_{k}}}
\newcommand\pin{\ifmmode \pi_n \else $\pi_{n}$\fi}
\newcommand\inv{\rm inv}

\begin{document}
\textheight=7.8truein
\runninghead{J.~A.~Kneissler}{Woven Braids and their Closures}
\normalsize\textlineskip
\thispagestyle{empty}
\setcounter{page}{1}

\title{\vspace{0.5cm}Woven Braids and their Closures\vspace{12pt}}
\author{\vspace{6pt}Jan A. Kneissler\\{\tenit Veilchenstrasse 6}\\{\tenit 75053 Gondelsheim}\\{\tenit Germany\vspace{16pt}}}
\maketitle
\vspace{0.5cm}

\abstracts{A special class of braids, called woven, is introduced and it is shown
qcxthat every conjugation class of the braid group contains woven braids.
In consequence tame links can be presented as plats and closures of woven braids.
Restricting on knots we get the `woven version' of the well-known theorem of Markov,
giving moves that are capable of producing all woven braids with equivalent closures.}
{As corollary we obtain that a link in which each component is dyed with at
least two different colors can be projected on a plane without crossing strands of
the same color.}
{Finally, a table of all minimal woven braids for the 84 prime knots with
at most nine crossings is appended. The average word length and the average number of entries per
knot type turn out to be surprisingly small.}

\vspace*{5pt}

%\textlineskip                  %) USE THIS MEASUREMENT WHEN THERE IS
%\vspace*{12pt}                 %) NO SECTION HEADING

\vspace*{1pt}\textlineskip

\section{Introduction and Results}\vspace*{-12pt}\subsection{Intention}
\setcounter{lemmac}{0}\setcounter{equation}{0}
%%\vspace*{-0.5pt}
%\firstline
\textheight=7.8truein
% \tableofcontents
% \newpage
%%\subsection{Intention}
\label{intention}
\firstline
The basic aim that led to these results was to look out for simple representatives
of conjugation classes of the braid group \Bn. In other words:~Find a family \Wn\ of braids
with $\forall\beta\in \Bn: \exists \gamma \in \Bn:
\gamma^{-1}\beta\gamma \in \Wn$! It is desirable to make this family $\Wn \subset \Bn$ in some sense as `small' as possible.

F.~Garside gave a very good solution to this problem in \cite{Ga}, constructing exactly one
representative per conjugation class, thereby solving the conjugacy problem in \Bn.
Unfortunately the characterization of these special braids
(called ``normal forms'') is not at all a trivial task (see theorem 2.7 of \cite {Bi}).

In this paper, looking for an other solution of the above-mentioned problem, we are
willing to sacrifice uniqueness for the sake of simplicity. The idea is to
try to let \Wn\ consist only of braids that can be presented with a maximum (compared to their conjugates)
number of constant strings (single straight lines).

\subsection{Woven braids}
\firstline
Given the usual generators $(\sigma_{i})_{i\in\{1,\ldots,n-1\}}$ of \Bn, letting the $i$th string pass over the $i\!+\!1$st, and $k$ positive integers $n_{1},\ldots,n_{k}$
with
%% $n = \sum \limits _{i=1}^{k} {n_{i}}$, we make the following abbreviations: (where $\hat{}$ means omission)
$n = n_1 + \cdots + n_k$, we make the following abbreviations: (\ $\hat{}$ means omission)
\begin{eqnarray}
%\nonumber
A_{ij}^{} & := & \sigma_{j-1}^{} \cdots \ \sigma_{i+1}^{} \sigma_{i}^{2} \sigma_{i+1}^{-1} \cdots \ \sigma_{j-1}^{-1} \qquad (1 \leq i < j \leq n)
\label{11}\\
%\nonumber
\pi_{n_{1} n_{2} \ldots n_{k}} & := & \sigma_1 \sigma_2 \cdots \ \hat \sigma_{n_{1}} \cdots \ \hat \sigma_{n_{1}+n_{2}} \cdots  \ \hat
 \sigma_{n_{1}+\cdots+n_{k-1}} \cdots \ \sigma_{n_{1}+\cdots+n_{k}-1}
\end{eqnarray}

\noindent
The $A_{ij}$ generate the pure braid group $\Pn := \langle (A_{ij})_{1 \leq i < j \leq n} \rangle$,
$\Pn^{j} := \langle (A_{ij})_{1 \leq i < j} \rangle$ is called {\it group of pure $j$-braids} and
\pink\ is a permutation braid.

From the theory of braids (e.g.\ \cite{Bi})
it is known that the $\Pn^{j}$ are free groups of rank $j-1$, that \Pn\ is the semidirect product of
$\Pn^{2},\ldots, \Pn^{n}$ and that every braid can be presented as product of a permutation braid
and a pure braid. We will now focus on special combinations of pure and permutation braids.

\skipline
\begin{definition}
\label{defwoven}
A braid on $n = n_1 + \cdots + n_k$ strings that can be presented in the form $\pink \beta_{n_{1}} \beta_{n_{1}+n_{2}}  \cdots \ \beta_{n_{1}+\cdots+n_{k}}$
with each $\beta_{i} \in \Pn^{i}$ is called a {\it woven} braid (with $k$ {\it components\/}) of {\it type} $(n_{1},\ldots,n_{k})$.
The set of all woven braids on $n$ strings with $k$ components will be called \Wnk.
\end{definition}
\skipline
Pictures of woven braids of type (4), (5) and (3,1,2) are given in
figure \ref{fig1}.\footnote{The second example explains our terminology, for this is how
typical outcomes of looms look like.}

\begin{center}
\fighead
\leavevmode
\showfbox{\epsfbox{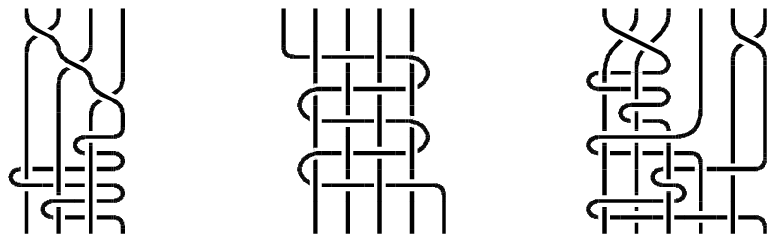}}
\figtail{Examples of woven braids.}
\label{fig1}
\end{center}

\noindent
We will prove in section 2 that $\Wn := \bigcup_{k=1}^{n}\Wnk$ indeed has the desired property:
\skipline
\begin{theorem}
\label{conj}
Every conjugation class of \Bn\ contains at least one element of\/ \Wn.
\end{theorem}

\subsection{Woven braids and links}
\firstline
For convenience, the closure of a
woven braid of the above type is named {\it web of index} $n$ and the plat
formed out of a woven braid whose type is a $k$-tuple of {\it even\/} integers with sum $2n$
is called {\it woven plat of index} $n$.\footnote{The closure of a braid is
obtained by identifying corresponding left and right endpoints; a plat is built by grouping
the endpoints on each side into non-separating pairs, and then glueing the members of each pair together.}
\pagebreak

\noindent
It is easily seen that both constructions yield a link with $k$ components.
We are interested only in tame links, and equivalence of link representatives shall
be given by ambient isotopy.
A link representative is called {\it multicolored} if each embedded circle is split up into
a finite collection of arcs, and a color is attached to each arc,
using at least two different colors per component.
With these notations we can state simple consequences of theorem \ref{conj}:
\skipline
\begin{corollary}
\label{corlink}
Let $L$ be an arbitrary equivalence class of links, then:
\begin{mylist}{}
\it
\item[{\rm $($a$)$}]
$L$ contains (infinitely many) webs and woven plats.
\item[{\rm $($b$)$}]
The minimal index of all webs in $L$ equals the braid index of $L$.
\item[{\rm $($c$)$}]
The minimal index of all woven plats in $L$ equals the bridge number of $L$.
\item[{\rm $($d$)$}]
For every multicolored representative $l$ of $L$ there exists a
color preserving ambient isotopy taking $l$ to $\tilde l$ and a projection of $\tilde l$ to a
colored regular planar diagram, in which the two strands of every crossing have different colors.
\end{mylist}
\end{corollary}

\skipline
Part (a) of this corollary is just a sharper version of Alexander's theorem (\cite{Al}), and one might ask for
a specialization of Markov's theorem (\cite{Ma} and \cite{Bi}) as well.
The following subsection gives the answer for the $k = 1$ case.

\subsection{Woven braids and knots}
\firstline
\begin{definition}
\label{defnabla}
Let $\nabla$ be the endomorphism of \Pn\ given by
\begin{equation}
%\mbox{Let $\nabla$ be the \Pn\ endomorphism given by}
\nabla (A_{ij}) := \left\{
\begin{array}{ll}
1_{\Pn} & $for $ i = 1 \\
A_{i-1,j-1} & $for $ i \geq 2\\
\end{array}
\right..
\label{13}
\end{equation}
To see that $\nabla$ is well-defined, check out that it is
compatible with the defining relations of \Pn\ that can be found in \cite{Mo}.
Furthermore let
\begin{equation}
\Kn := \lbrace \ \beta \, \nabla(\beta) \, \nabla\big(\nabla(\beta)\big) \cdots \nabla^{n-2}(\beta)\  \vert \ \beta \in \Pn^{n} \  \rbrace
\end{equation}
and define the following moves for woven braids:

{\vspace{6pt}
\begin{tabular}{p{2.4cm}l}
Type{ \hfill I\hspace{6pt}\hfill}move: &
Replace $\omega \in \Wne$ by $\kappa^{-1}\omega\kappa$ for some $\kappa \in \Kn$.
\vspace{1pt}\\
Type{ \hfill II$^+$\hfill}move: &
Replace $\omega \in \Wne$ by $\omega\sigma_{n}^{\pm 1}$.
\vspace{1pt}\\
Type{ \hfill II$^-$\hfill}move: &
Replace $\omega\sigma_{n-1}^{\pm1} \in \Wne$ by $\omega$ when $\omega \in {\cal W}_{n-1}^1$. \\
\end{tabular}
\vspace{6pt}}
\end{definition}

\skipsline
\remark
An explicit specification of the string number $n$ is not necessary here, for the
permutation of a braid in \Wne\
%%% , which can be read off the braid word,
is a cycle of length $n$.

\skipline
\begin{theorem}\label{knots}
\it
\begin{mylist}{}
\item[{\rm $($a$)$}]
Application of a type I and II$^{\pm}$ move on a woven braid of type $(n)$ yields a woven braid of type $(n)$ and $(n\pm 1)$, respectively.
\item[{\rm $($b$)$}]
\Kn\ is a subgroup of \Pn\, operating on \Wne\ by conjugation; the orbits are
the intersections of the conjugation classes of \Bn with \Wne.
\item[{\rm $($c$)$}]
Two webs represent the same knot, if and only if the corresponding woven braids are related via a
finite sequence of type I, II$^{\pm}$ moves.
\end{mylist}
\end{theorem}

\section{Weaving Braids by Conjugation}
\setcounter{lemmac}{0}\setcounter{equation}{0}
\label{proof1}
\firstline
If for some fixed $s \in \{2,\ldots,n\}$
the pure braid $\beta$ can be written as word of \Pn\ without using letters $A_{is}$  i.e.\ $\beta \in \langle(A_{ij})_{j\neq s}\rangle$,
then we say string $s$ is {\it free} in $\beta$.

\skipline
\begin{lemma}
\label{free}
For any $\beta \in \Pn$ and $s \in \{2,\ldots,n\}$ there exist two unique pure $s$-braids $_s\beta,\beta_s \in \Pn^{s}$
such that string $s$ is free in $(_s\beta)^{-1}\beta$ and in $\beta(\beta_s)^{-1}$.
\end{lemma}
\skipsline
\remark
$\beta_s = {_s\beta} = 1_{\Pn}$ if and only if string $s$ is free in $\beta$.
\skipsline
\proof{}
Due to the structure of \Pn, being the semidirect product of the $\Pn^k$, pure $i$-braids commute
with pure $j$-braids in the sense that for all $i < j$, $\mu \in \Pn^{i}, \nu \in \Pn^{j}$
there exists a $\tilde\nu \in \Pn^{j}$ such that $\nu\mu = \mu\tilde\nu$.
Thus one can push all subwords $\in \Pn^{s}$ of a word for $\beta$ to the beginning
or to the end where they can be cancelled.

\noindent
Uniqueness: If we have two such $\grave\beta_s, \acute\beta_s \in \Pn^{s}$ then
string $s$ has to be free in
$\big(\beta{\grave\beta_s}^{-1}\big)^{-1} \beta{\acute\beta_s}^{-1} = \grave\beta_s^{}{\acute\beta_s}^{-1}$
as well. This implies $\grave\beta_s = \acute\beta_s$ (analogously for $_s\beta$)\qedp

\skipline
\begin{lemma}
\label{comb}
For any $\beta \in \Pn$: $\beta = \beta_2\ldots\beta_n = {_n\beta}\ldots{_2\beta}$.
\squad{\rm (``Combing'')}
\end{lemma}

\skipsline
\proof{}
Arrange the subwords $\in \Pn^{j}$ with ascending $j$ by commuting them
as described in lemma 2.\ref{free}. The result is $\check\beta_2\ldots\check\beta_n$
with all $\check\beta_j \in \Pn^{j}$. Pushing $\check\beta_s$ to the end
leaves it unchanged, because it has to pass only through those $\check\beta_j$ with $j > s$,
thus $\check\beta_s = \beta_s$. The statement for descending order is proved correspondingly\qedp

\skipline
\proof{ of theorem \ref{conj}}
\skipsline
Given an arbitrary element $\beta \in \Bn$, we are looking for a conjugate $\hat\beta \in \Wn$
of $\beta$.\\
Let $\Pi : \Bn \to S_n$ denote the homomorphism given by $\Pi(\sigma_i) = (i \quad i+1)$,
then $\Pn = \ker\Pi$.
Say $\Pi(\beta)$ consists of $k$ cycles of the lengths $n_1,\ldots,n_k \geq 1$
with $n = n_1+\cdots+n_k$.
Since the cycle type determines conjugation classes in $S_n$, there exists a permutation
$\tau \in S_n$ with
\begin{eqnarray}
\tau^{-1}\Pi(\beta)\tau & = &
(n_1 \; n_1\!-\!1 \; \ldots \; 1)(n_1\!+\!n_2 \; \ldots \; n_1\!+\!1)\cdots
({\textstyle \sum_{i=1}^{k}\limits} n_i \; \ldots \; {\textstyle \sum_{i=1}^{k-1}\limits} n_i\!+\!1)
\nonumber\\
 & = &\Pi(\pink).
\label{pireq}
\end{eqnarray}
\noindent
Take any $\pi_\tau \in \Bn$ with $\Pi(\pi_\tau) = \tau$ and let
$\tilde\beta := \pink^{-1}\pi_\tau^{-1}\beta\,\pi_\tau^{}$,
then $\Pi(\tilde\beta) = 1_{S_n}$ thus ~$\tilde\beta \in \Pn$.
Now let
\begin{eqnarray}
\Delta (\alpha) & := & \mathpink^{}\alpha\,\mathpink^{-1} \label{26o} \\
N & := & \{n_1, n_1\!+\!n_2, \ldots, n_1\!+\cdots+\!n_k\}, \quad \bar N := \{1,\ldots,n\} \setminus N \\
F(\alpha) & := & \big\{\: s \; \vert \; \mbox{string }s\mbox{ is not free in }\alpha^{-1}\tilde\beta\:\Delta (\alpha)\, \big\} \cap \bar N
\label{26u}
\end{eqnarray}
then $\Delta$ is an automorphism of \Pn\ and $F$ is a set-valued function defined on $\Pn$
because for all $\alpha \in \Pn$ we have $\Pi\big(\alpha^{-1}\tilde\beta\:\Delta(\alpha)\big) = \Pi(\tilde\beta) = 1_{S_n}$.
\pagebreak

\noindent
We will now prove the existence of an $\alpha \in \Pn$ with $F(\alpha) = \emptyset$ inductively:
Supposing there is an $\alpha_{(s-1)}$ with $\inf F\big(\alpha_{(s-1)}\big) > s-1$,
we construct $\alpha_{(s)}$ with $\inf F\big(\alpha_{(s)}\big) > s$.
Starting with $\alpha_{(0)} := 1_{\Pn}$, we thereby end up with $\alpha := \alpha_{(n)}$, having the
property $\inf F\big(\alpha_{(n)}\big) > n$; this implies $F(\alpha) = \emptyset$
\quad ($\inf \emptyset := +\infty$).
%.\footnote{$\inf \emptyset = +\infty$.}
\vspace{2pt}\begin{mylist}{\setlength{\leftmargin}{15pt}}
\item[1.]
$\inf F\big(\alpha_{(s-1)}\big) > s$: Simply take $\alpha_{(s)} := \alpha_{(s-1)}$.
\vspace{2pt}
\item[2.]
$\inf F\big(\alpha_{(s-1)}\big) = s$: This can only happen for $s \in \bar N$. Now let
\begin{eqnarray}
\delta & := & (\alpha_{(s-1)})^{-1}\,\tilde\beta\:\Delta(\alpha_{(s-1)})\\
\alpha_{(s)} & := & \alpha_{(s-1)}\:{(_s\delta)}\label{eq27}\\
\delta^\prime & := & (\alpha_{(s)})^{-1}\,\tilde\beta\:\Delta(\alpha_{(s)}) = ({_s\delta})^{-1}\,\delta\,\Delta({_s\delta}).
\end{eqnarray}
String $s$ is free in $({_s\delta})^{-1}\delta$ (lemma 2.\ref{free}).
Furthermore the strings with numbers in $\{1,\ldots,s-1\} \cap \bar N$ are free in both $\delta$ (by
induction hypothesis) and in $_s\delta$ (because $_s\delta \in \Pn^{s}$), hence they are free in $({_s\delta})^{-1}\,\delta$, too.

For $s \in \bar N$ we have $\Delta (\Pn^{s}) \subset \Pn^{s+1}$, which can easily be verified,
looking at the generators $A_{is}$; the three typical cases are shown in figure \ref{fig2}a--c.
So the strings with numbers $\leq s$ are also free in $\Delta({_s\delta})$.
Thus all strings with numbers in $\{1,\ldots,s\} \cap \bar N$ are free in $\delta^\prime$ and
(\ref{26u}) yields $\inf F\big(\alpha_{(s)}\big) > s$.
\end{mylist}
\skipline
Now, having $F(\alpha) = \emptyset$ and $\tilde\beta$ as above, we can construct $\hat\beta$ as follows:
\begin{equation}
\hat\beta := \pink\alpha^{-1}\tilde\beta\,\Delta(\alpha)
= \big(\pi_\tau\,\Delta(\alpha)\big)^{-1}\beta\,\big(\pi_\tau\,\Delta(\alpha)\big).
\label{eq28}
\end{equation}
With respect to the first expression, lemma 2.\ref{comb} and
the remark to lemma 2.\ref{free} together with $F(\alpha) = \emptyset$
show that $\hat\beta$ is of the
form $\pink \: \hat\beta_{n_{1}} \, \hat\beta_{n_{1}+n_{2}}  \cdots \ \hat\beta_{n_{1}+\cdots+n_{k}}$.
The last term of (\ref{eq28}) makes clear that $\hat\beta$ is conjugate to $\beta$\qedp
\begin{center}
\fighead
\leavevmode
\showfbox{\epsfbox{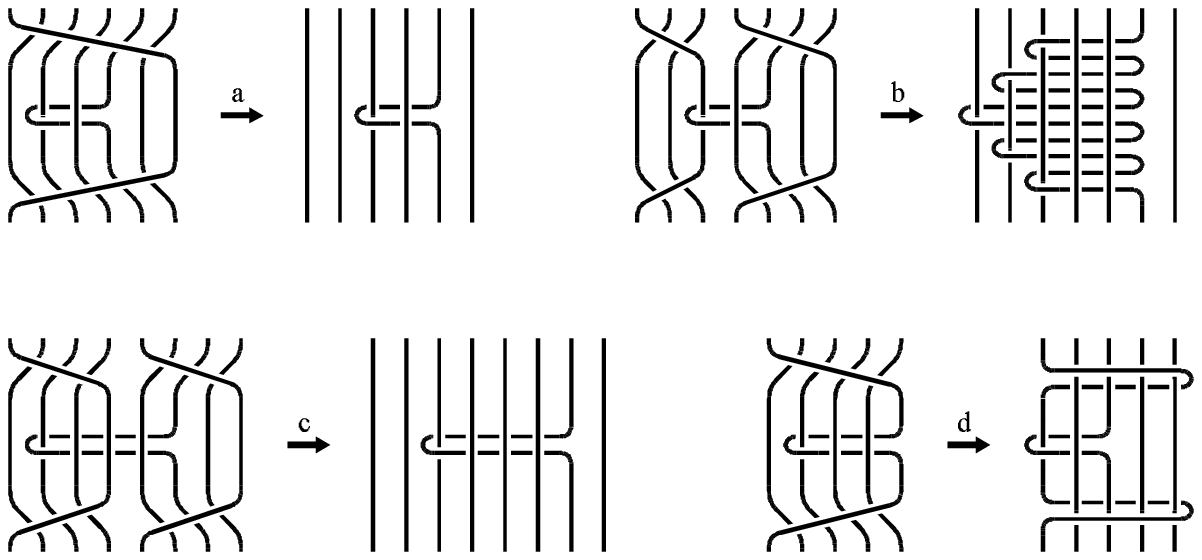}}
\vspace*{-0.5cm}
\figtail{The action of $\Delta$.}
\label{fig2}
\end{center}
\pagebreak
\proof{ of corollary \ref{corlink}}
\skipsline
$($a$)$ and $($b$)$: Besides the results of \cite{Al} and \cite{Ma} there is only to mention
that closing a woven braid of type $(n_1,\ldots,n_k)$,
placed parallel to the $x$-$y$-plane,
around the $x$-axis and projecting it along the $z$-axis, yields a woven plat of
type $(2n_1,\ldots,2n_k)$. Hence every web of index $n$ can be regarded
as woven plat of index $n$ (but not vice versa).
\skipsline
$($c$)$: If the bridge number of $L$ is $m$, there exists a $2m$-plat representing
$L$. We have to transform this plat into a
woven plat without changing the number of strings. This is a task quite similar
to the proof of theorem \ref{conj}, so we will give only the crucial points here.
Again, the first step is to order the strings by adjusting the permutation of the braid to
$\Pi(\pink)$. This can be done by using moves of type a and b in
figure \ref{fig3} at {\it both} ends of the plat.

In order to get rid of unwanted $A_{ij}$'s, we pushed them
upwards and all around the closed braid. Passing $\pink$ increased the index $j$,
with the result that the $A$'s assembled at the highest available string of
each component (those with numbers in $N$).
In the case of plats we have to push $A_{ij}$'s
with even $j \in \bar N$ upwards, they thereby have to pass \pink\ and arrive as
elements of $\Pn^{j+1}$
at the upper end. The $A_{ij}$'s with odd $j$ have to be pushed down to the lower end of the plat.
In both cases we can then apply the moves c and d of figure \ref{fig3} or their mirror images
repeatedly, in order to free the strings with numbers in $\bar N$ successively, just as before.
\begin{center}
\fighead
\leavevmode
\showfbox{\epsfbox{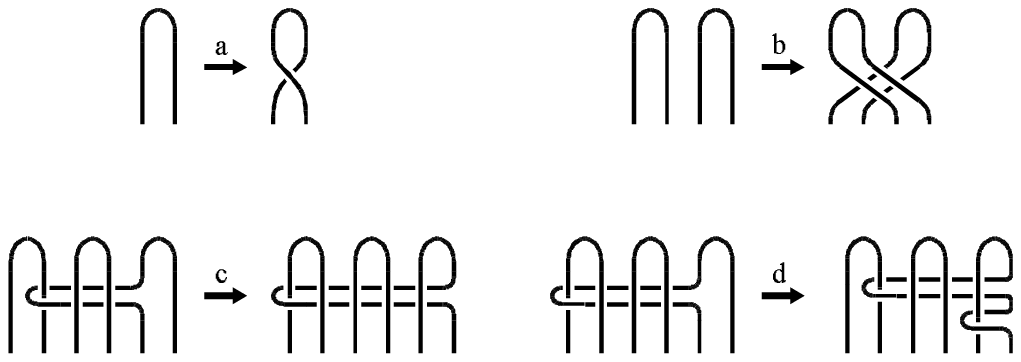}}\figtail{Useful moves for plats, preserving the link type.}\label{fig3}
\end{center}
\noindent
\begin{minipage}[t]{3.7truein}
\skipsline
$($d$)$: Having more than two arcs per component or more than two colors makes finding an
appropriate projection only easier, so we may assume that each embedded circle
consists of exactly two arcs, colored black and gray.
Choose $\tilde l$ to
be a web equivalent to $l$, presented as required in definition \ref{defwoven}. By parameter transformation we may arrange
the coloring of $\tilde l$ in the way indicated in figure \ref{fig4}, identifying the black arcs
with the `moving' parts of the strings with numbers $\in N$. The obtained
regular diagram has no monochrome crossings\qedp
\end{minipage}
\hspace{\fill}
\begin{minipage}[t]{1.1truein}\vspace{9pt}\begin{center}\leavevmode\showfbox{\epsfbox{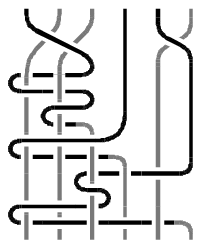}}
\vspace*{7pt}\fcaption{}\label{fig4}\end{center}\end{minipage}
\vspace{16pt}
\pagebreak
\section{The algebraic Knot Problem for Webs}
\setcounter{lemmac}{0}\setcounter{equation}{0}
\firstline
For a braid $\beta \in \Bn$ placed in the $x$-$y$-plane we will denote by $\inv(\beta)$
the braid obtained by $\pi$-rotating $\beta$ on the $z$-axis. Observe that
this rotation reverses the order of the braid word and sends $\sigma_i$
to $\sigma_{n+1-i}$, so inv is an involution and an antihomomorphism.

\skipline
\begin{lemma}
\label{inv}
If $\omega \in \Wne$ then $\inv(\omega) \in \Wne$.
\end{lemma}
\skipsline
\proof{} Showing this algebraically is an easy task, so we will content ourselves
with the following argumentation: The elements of \Wne\ are characterized by the fact
that they can be presented as braids with $n-1$ constant strings and one
string moving arbitrarily from position 1 to position $n$; a property that is preserved under
inv\qedp
\skipsline
\remark The analogous statement for $\Wnk$ is false when $k \geq 2$.
\skipline
We will assume from now on that $k = 1$.
Composing the $\Pn$-endomorphism $\nabla$ of definition \ref{defnabla} and
the $\Pn$-automorphism $\Delta (\alpha) = \pin^{}\alpha\,\pin^{-1}$ of eq.~(\ref{26o})
defines an endomorphism $\Phi_n := \nabla\circ\Delta\,:\, \Pn \to \Pn$ that
will be useful for the characterization of \Wne\ and \Kn.
\skipline
\begin{lemma}
\label{char}
$\Wne = \pin \ker(\Phi_n)$ and
$\Kn = \{ \kappa \in \Pn \; \vert \; \nabla(\kappa) = \Phi_n(\kappa) \}$.
\end{lemma}
\skipsline
\proof{} From figure \ref{fig2}a,d we obtain (with an appropriate $\chi \in \Pn$)
\begin{displaymath}
\Delta (A_{ij}) = \left\{
\begin{array}{ll}
A_{i+1,j+1} &\mbox{for }j < n \\
\chi^{-1}A_{1i}\:\chi &\mbox{for }j = n
\end{array}
\right.
\mbox{ and thus \squad}\Phi_n (A_{ij}) = \left\{
\begin{array}{ll}
A_{ij} &\mbox{for }j < n \\
1_{\Pn} &\mbox{for }j = n
\end{array}
\right..
\end{displaymath}
So $\ker(\Phi_n) = \Pn^{n}$, justifying the first statement.
For any $\beta \in \Pn^{n}$ we get:\footnote{According to eq.~(\ref{13}) we have $\nabla (\Pn^{j}) = \Pn^{j-1}$ for $j \geq 3$ and
$\nabla^{n-1}(\beta) = 1_{\Pn}$.}
\begin{equation}
\Phi_n \big(\,\beta\: \nabla(\beta) \cdots \nabla^{n-2}(\beta)\big)
= \nabla(\beta) \cdots \nabla^{n-2}(\beta) = \nabla \big(\,\beta\: \nabla(\beta) \cdots {\nabla}^{n-2}(\beta)\big).
\end{equation}
\noindent
Conversely, by combing $\kappa \in \Pn$ in descending order $\kappa = {_n\!\kappa}\cdots{_2\kappa}$
we get
\begin{eqnarray}
\nabla(\kappa) = \Phi_n(\kappa) &
 \Rightarrow &
\nabla({_n\!\kappa})\cdots\nabla({_3\kappa}) = {_{n-1}\!\kappa}\cdots{_2\kappa}.
\quad \Rightarrow \quad
{_{i-1}\!\kappa} = \nabla({_{i}}\kappa) \nonumber\\
& \Rightarrow & \kappa = {_n\!\kappa}\;\nabla({_n\!\kappa})\;\nabla\big(\nabla({_n\!\kappa})\big)\cdots\nabla^{n-2}(_n\!\kappa) \in \Kn.
\end{eqnarray}
\noindent
In the second equation both sides are descendingly combed, so the second
implication is just the uniqueness condition of lemma 2.\ref{free} for $2 < i \leq n$\qedp
\skipline
\proof{ of theorem \ref{knots}}
\skipsline
The simple relation $\sigma_n^{}A_{in}^{}\sigma_n^{-1} = A_{i,n+1}^{}$ in ${\cal B}_{n+1}^{}$
shows for $\beta \in \Pn^{n}$ that
\begin{equation}
\pin^{}\beta\,\sigma_n^{-1} = \pi_{n+1}^{}\,\sigma_n^{-2}\,\beta^\uparrow\mbox{\quad and\quad}
\pin^{}\beta\,\sigma_n^{} = \pi_{n+1}^{}\,\sigma_n^{-2}\,\beta^\uparrow\sigma_n^2
\label{312}
\end{equation}
where $\beta^\uparrow \in {\cal P}_{n+1}^{n+1}$ is obtained from $\beta$ by augmenting the second index of all
$A_{ij}$'s by one. Since $\sigma_n^2 = A_{n,n+1}^{}$, we see that a type II$^+$ move
yields a woven braid with incremented index.

\pagebreak
\noindent
We have for $\omega \in \Wne$ and $\kappa \in \Pn$ the following consequence of lemma 3.\ref{char}:
\begin{eqnarray}
{\kappa^{-1}\omega\kappa \in \Wne} & {\Leftrightarrow} & {\pin^{-1}\kappa^{-1}\omega\kappa = \Delta^{-1}(\kappa^{-1})\pin^{-1}\omega\kappa \in \ker \Phi_n}
\nonumber\\
{} & {\Leftrightarrow} & {\Phi_n\big(\Delta^{-1}(\kappa^{-1})\big)\Phi_n(\pin^{-1}\omega)\Phi_n(\kappa) =
\nabla(\kappa^{-1})\Phi_n(\kappa) = 1_{\Pn}}
\nonumber\\
& \Leftrightarrow & \nabla(\kappa) = \Phi_n(\kappa) \quad \Leftrightarrow \quad \kappa \in \Kn
\label{313}
\end{eqnarray}
The ``$\Leftarrow$''-direction states that the result of a type I move is again a woven braid.
Conversely we get that \Kn\ is a subgroup of \Pn\, because for any
$\kappa_1,\kappa_2 \in \Kn$ the woven braids $\kappa_1^{-1}\omega\kappa_1$
and $\kappa_2^{-1}\omega\kappa_2$ are conjugate by $\kappa_1^{-1}\kappa_2 \in \Pn$
and the ``$\Rightarrow$''-direction thus implies $\kappa_1^{-1}\kappa_2 \in \Kn$.
\skipsline
Let $\omega, \gamma^{-1}\omega\gamma \in \Wne$ with $\gamma \in \Bn$.
Then $\Pi(\gamma^{-1}\omega\gamma) = \Pi (\omega) \Rightarrow
\big[\Pi (\gamma),\Pi (\omega)\big] = 1_{S_n}$; but since $\Pi(\omega) = (n\; \ldots \; 1)$ consists
of exactly one cycle, this is only possible if $\Pi(\gamma) = \big(\Pi (\omega)\big)^i$
with some integer $i$.
Thus we have a pure braid $\tilde\gamma := \omega^{-i}\gamma \in \Pn$ giving the
same conjugation result: $\gamma^{-1}\omega\gamma = \tilde\gamma^{-1}\omega\tilde\gamma$.
According to (\ref{313}) $\tilde\gamma \in \Kn$, showing
that any two conjugate woven braids are related via a single type I move.
\skipsline
Markov's theorem states that the closures of two braids represent the same link if and only
if one braid can be transformed into the other by conjugation and stabilization-moves.
So it only remains to show, that for any $\beta \in \Bn$ and
$\beta^\prime = \beta\sigma_n^{\pm1}\in {\cal B}_{n+1}$ there are woven
braids $\omega \in \Wne$ and $\omega^\prime \in {\cal W}_{n+1}^1$ that are related
by a type II$^\pm$ move and are
contained in the conjugation classes of $\beta$ and $\beta^\prime$, respectively.

Similar to the proof of theorem \ref{conj} we construct a $\gamma \in \Bn$
with $\gamma^{-1}\inv(\beta)\gamma \in {\cal W}_{\it n}^1$.
According to equations (\ref{eq28}) and (\ref{eq27}) we obtain $\gamma = \pi_\tau\Delta(\alpha)$ where
$\alpha$ contains only the generators $A_{ij}$ with $j < n$.
Therefore the action of $\Delta$ on $\alpha$ is described by figure \ref{fig2}a, and
thus $\Delta(\alpha) \in \langle (A_{ij})_{2\leq i<j\leq n}\rangle$.

The requirement $\tau^{-1}\Pi(\inv(\beta))\tau = \Pi(\pi_{\it n})$ of equation (\ref{pireq}) determines
the permutation $\tau$ not uniquely: any other $\tilde\tau := (\Pi(\pin))^i\,\tau$
works as well. Since $\Pi(\pin)$ is a cycle of length $n$, there is an integer $i$
such that $\tilde\tau (1) = 1$ and we can choose $\pi_\tau \in \langle (\sigma_i)_{2 \leq i < n}\rangle$.
\skipsline
We have therefore found a $\gamma$ that can be written as word of \Bn\ without
using $\sigma_1$. So in $\tilde\gamma := \inv(\gamma)$ no $\sigma_{n-1}$'s
appear and we have  $[\tilde\gamma,\sigma_n] = 1_{{\cal B}_{n+1}}$
in ${\cal B}_{n+1}$.
We can now take
\begin{eqnarray}
\omega & := & \inv (\gamma^{-1}\inv(\beta)\gamma) = \tilde\gamma\beta\tilde\gamma^{-1}
\mbox{\quad and\quad} \nonumber\\
\omega^\prime & := & \tilde\gamma\beta^\prime\tilde\gamma^{-1} =
\tilde\gamma\beta\sigma_n^{\pm1}\tilde\gamma^{-1} =
\tilde\gamma\beta\tilde\gamma^{-1}\sigma_n^{\pm1} = \omega\sigma_n^{\pm1}.
\end{eqnarray}
By lemma 3.\ref{inv} $\omega$ is a woven braid and $\omega^\prime$ is obviously the result
of a type II$^+$ move applied on $\omega$. This completes the proof of theorem \ref{knots}\qedp
\skipline
\remark
The premise $k = 1$ was very helpful at several points of this proof and
things are getting much more complicated for links. For $k > 2$ conjugations with
non-pure braids are inevitable because one must have moves that exchange
components.
It is likely, however, that -- with some additional moves -- theorem \ref{knots} can be
extended to the case of colored links.
%%%%% , where the components carry different labels.

\section{Tables}
\setcounter{lemmac}{0}\setcounter{equation}{0}
\vspace*{-12pt}
\subsection{Link-tabulation with woven braids}
\firstline Due to corollary \ref{corlink} links can be tabulated by giving a woven braid
$\omega_L$ for each equivalence class $L$, with either the closure of $\omega_L$ or the
plat formed out of $\omega_L$ is contained in $L$. Depending on the construction
we will call $(\omega_L)_{L \in {\cal L}}$ a web- or a plat-table for the family
$\cal L$ of links.
We have the nice side effect that orientation and
labelling of the components are provided automatically, if we agree upon the convention that all
strings of the closure and all odd-numbered strings of the plat shall be oriented in reading direction
of the braid word and that components shall be labelled increasingly with respect to their position.

In comparison with tables based on diagrams we lose the information about the crossing numbers and
alternation type.\footnote{For instance, one can show that any web representing
the figure-eight knot must have at least six crossings and cannot be alternating.}
 On the other side, web-/plat-tables can implicitly convey the braid-/bridge-indices and in some
cases they reveal strong symmetries. A knot is called strongly invertible /
strongly negative amphicheiral if there is an orientation preserving / reversing
self-homeomorphism of $S^3$ that is an involution fixing the knot
setwise and reversing the knot orientation.
A woven braid $\omega$ is called {\it symmetric} if $\mbox{inv}(\omega) = \omega$ and
{\it antisymmetric} if $\mbox{inv}(\omega) = \mbox{mirr}(\omega)$, with the automorphism mirr given by $\mbox{mirr}(\sigma_i^{}) := \sigma_i^{-1}$.
It is a simple consequence that a web or woven plat build out of a symmetric / antisymmetric
woven braid is strongly invertible / strongly negative amphicheiral.

%\footnote{Webs having the second symmetry appeared in table \ref{knottable}
%automatically for all prime amphicheiral
%knots with at most nine crossings.}

\subsection{Notation}
\firstline
For woven braids we have a trivial solution of the word problem of \Bn:
Let the woven braid be presented in the form required in definition \ref{defwoven},
expand all $A_{ij}$ according to (\ref{11}) and then reduce until no
more subwords of type $\sigma_i^{\pm1}\sigma_i^{\mp1}$ occur.
The result, which will be called {\it tight} word, is unique because
the $\Pn^{i}$ are free. In general, tight words are not
minimal in length; $(\sigma_2^{}\sigma_1^{-1})^2$ for example yields the same braid as the
tight word $\sigma_1^{-1}\sigma_2^{-2}\sigma_1^2\sigma_2^{}$.

Tight words could in principle be written down by listing the indices and exponents
of the $\sigma_i$'s, but we will use a much more efficient notation for knots making use
of the fact that in tight words of $\Wne$ every $\sigma_i$ involves the string starting at position 1.
\skipline
\begin{lemma}
For any tight word $w = \sigma_{i_1}^{e_1}\ldots\sigma_{i_l}^{e_l}$ of \Wne\ with
$1 \leq i_t < n$ and $e_t = \pm1$ the $l$-tuple $(e_1,c_1,c_2,\ldots,c_{l-1})$ determines $w$ uniquely, if we take
\begin{equation}
c_t := \left\{
\begin{array}{ll}
0 & $if $ i_t \neq i_{t+1} $ and\/ $ e_t = e_{t+1} \\
1 & $if $ i_t \neq i_{t+1} $ and\/ $ e_t \neq e_{t+1} \\
2 & $if $ i_t = i_{t+1}
\end{array}
\right\}
\quad \mbox {for\/ } 1 \leq t \leq l-1.
\end{equation}
\end{lemma}
%\skipsline
%\pagebreak

\proof{}
The indices and exponents can be recovered recursively:
\begin{equation}
\begin{array}{rclrcl}
d_1 & := & 1, & d_{t+1} & := & (-1)^{\frac{1}{2}c_t(c_t-1)}d_t \\
i_1 & := & 1, & i_{t+1} & := & i_t + \big(1 - \frac{1}{2}c_t(c_t-1)\big)d_t \\
    &    &    & e_{t+1} & := & (-1)^{c_t}e_t.
\end{array}
\end{equation}
This statement is obviously correct for the $e_t$.
We call the string moving from first position in the braid to position $n$ the {\it main} string.
Observing that $\frac{1}{2}c_t(c_t-1) = 1$ for $c_t = 2$ and $= 0$ otherwise, we realize that $d_t$ describes
the direction of the main string at the crossing $\sigma_{i_t}^{e_t}$. Indeed $d_1 = 1$
(towards higher indices) at the beginning of the braid, and a sign change takes place
whenever two consecutive $\sigma$'s are at the same level (at U-turns).
The index $i_t$ is increased or decreased whenever $c_t \neq 2$, depending on the
current direction $d_t$\qedp
\skipline
In web-tables we can even strip off leading and tailing zeros and use
the tuple $(e_1,c_u,\ldots,c_v)$
instead, where $u$ and $v$ are min and max of $\{j\; \vert\; c_j \neq 0\}$, respectively.
The reason is that if the tuples of woven braids differ only by the numbers
of their leading and tailing zeros, then the braids are related by type
II$^\pm$ and $\inv \circ \mbox{II}^\pm \circ \inv$ moves and their closures are of the same knot
type. The tuple $(e_1,c_u,\ldots,c_v)$ shall therefore refer to the one with lowest string index,
which is uniquely determined.
\skipline
\remark
\begin{mylist}{}
\item
We have hereby reduced braid words of length $l$ with $2(n-1)$ different letters
to words of length $\leq l$ using only three letters.
\item
If $(e_1,c_u,\ldots,c_v)$ is the tuple describing $\omega$  then
$\mbox{mirr}(\omega)$ is given by the tuple $(-e_1,c_u,\ldots,c_v)$ and $\inv(\omega)$ by
$((-1)^{c_u+\cdots+c_v}e_1,c_v,\ldots,c_u)$.
\item
An analogous notation for links can be given as well by adding the
letter ``3'', signalling the passage to the next component.
\end{mylist}

\skipline
\subsection{Minimal webs}
\vspace*{-6pt}
\noindent
\begin {definition}
We call the closure of a woven braid {\it minimal}, if and only if
\begin{mylist}{}
\item[1.]
the number of strings equals its braid index, and
\item[2.]
its length (i.e. the length of its tight braid word) is minimal amongst all other webs
satisfying condition 1 and representing the same oriented (colored) link.
\end{mylist}
\end{definition}
\skipsline
\remark Minimal webs have not necessarily the minimal length of all equivalent
webs: Computer calculations revealed that
$\sigma_1^{}\sigma_2^{-2}\sigma_1^2\sigma_2^{}\sigma_3^{-2}\sigma_2^{}\sigma_1^{-2}\sigma_2^{-2}\sigma_1^2\sigma_2^{}\sigma_3^{}
 \in {\cal W}_4^1$ gives a minimal
web representing the knot $9_{34}$, but closing the shorter woven braid
$\sigma_1^{}\sigma_2^{}\sigma_3^{2}\sigma_2^{}\sigma_1^{-2}\sigma_2^{}\sigma_3^{-1}\sigma_4^{-2}\sigma_3^{}\sigma_2^{-2}\sigma_3^{}\sigma_4^{}
\in {\cal W}_5^1$
yields $9_{34}$ as well.
\skipline
Table \ref{knottable} lists {it all} minimal webs for the prime knots with at most nine
crossings.
The given woven braids are Markov-equivalent to the braids
listed in \cite{Jo}.\footnote{This condition determines the orientations and which knot out of chiral pairs is chosen.
Take into account that in \cite{Jo} the mirror images of the $\sigma_i$'s are taken
as generators of \Bn!}

%\begin{table}
{\footnotesize
\tcaption{Minimal webs for prime knots with crossing number $\leq 9$.}
\label{knottable}
%\begin{center}
\vspace{13pt}
\noindent
\mbox{
\begin{tabular}{|r@{}l|l|r@{\hspace{1ex}}r@{\hspace{0,5ex}}r|c|} \hline
\multicolumn{2}{|c|}{AB}&\multicolumn{1}{c|}{Webs}&\multicolumn{1}{c|}{b}&\multicolumn{1}{c|}{l}&\multicolumn{1}{c|}{n}&s\\
\hline
$3$&$_{1}$&${\it-5}^{\rm s}$&$2$&$3$&$1$&r\\
$4$&$_{1}$&$\pm15^{\rm a}$&$3$&$6$&$2$&f\\
$5$&$_{1}$&${\it-53}^{\rm s}$&$2$&$5$&$1$&r\\
$5$&$_{2}$&$-31$&$3$&$6$&$2$&r\\
$6$&$_{1}$&$39$&$4$&$9$&$2$&r\\
$6$&$_{2}$&$-143,-419$&$3$&$8$&$4$&r\\
$6$&$_{3}$&$\pm{\it100}^{\rm a}$&$3$&$6$&$2$&f\\
$7$&$_{1}$&${\it-485}^{\rm s}$&$2$&$7$&$1$&r\\
$7$&$_{2}$&$-91$&$4$&$9$&$2$&r\\
$7$&$_{3}$&$-319$&$3$&$8$&$2$&r\\
$7$&$_{4}$&$1794^{\rm s}$&$4$&$9$&$1$&r\\
$7$&$_{5}$&$-287,-851$&$3$&$8$&$4$&r\\
$7$&$_{6}$&$-99$&$4$&$9$&$2$&r\\
$7$&$_{7}$&$-3741^{\rm s},16986^{\rm s}$&$4$&$11$&$2$&r\\
$8$&$_{1}$&$111$&$5$&$12$&$2$&r\\
$8$&$_{2}$&$-1295,-3779$&$3$&$10$&$4$&r\\
$8$&$_{3}$&$\pm115^{\rm a}$&$5$&$12$&$2$&f\\
$8$&$_{4}$&$-359,-1067$&$4$&$11$&$4$&r\\
$8$&$_{5}$&$-3875$&$3$&$10$&$2$&r\\
$8$&$_{6}$&$-395,-1175,$&&&&\\
&&$-3515$&$4$&$11$&$6$&r\\
$8$&$_{7}$&${\it916}$&$3$&$8$&$2$&r\\
$8$&$_{8}$&$-280$&$4$&$9$&$2$&r\\
$8$&$_{9}$&$\pm1439^{\rm a},\pm4031$&$3$&$10$&$6$&f\\
$8$&$_{10}$&${\it-911}$&$3$&$8$&$2$&r\\
$8$&$_{11}$&$375$&$4$&$11$&$2$&r\\
$8$&$_{12}$&$\pm135^{\rm a}$&$5$&$12$&$2$&f\\
$8$&$_{13}$&$1814$&$4$&$9$&$2$&r\\
$8$&$_{14}$&$7629,-16970$&$4$&$11$&$4$&r\\
$8$&$_{15}$&$-2799,-3701,$&&&&\\
&&$-7699$&$4$&$11$&$6$&r\\
$8$&$_{16}$&$-12521^{\rm s}$&$3$&$10$&$1$&r\\
$8$&$_{17}$&$1275^{\rm an},-2715^{\rm n},$&&&&\\
&&$-3772^{\rm n},12083^{\rm an}$&$3$&$10$&$4$&i\\
$8$&$_{18}$&$\pm7708^{\rm a}$&$3$&$10$&$2$&f\\
$8$&$_{19}$&$-365$&$3$&$8$&$2$&r\\
$8$&$_{20}$&$371$&$3$&$8$&$2$&r\\
$8$&$_{21}$&$-905$&$3$&$8$&$2$&r\\
$9$&$_{1}$&${\it-4373}^{\rm s}$&$2$&$9$&$1$&r\\
$9$&$_{2}$&$-271$&$5$&$12$&$2$&r\\
$9$&$_{3}$&$-2911$&$3$&$10$&$2$&r\\
$9$&$_{4}$&$-955$&$4$&$11$&$2$&r\\
$9$&$_{5}$&$14106$&$5$&$12$&$2$&r\\
$9$&$_{6}$&$-2591,-7667$&$3$&$10$&$4$&r\\
$9$&$_{7}$&$-827,-2471,$&&&&\\
&&$-7403$&$4$&$11$&$6$&r\\
$9$&$_{8}$&$-291$&$5$&$12$&$2$&r\\
\hline
\end{tabular}
%\end{center}
%\end{table}
%\begin{table}
%\setcounter{table}{0}
%\tcaption{Continued.}
%\begin{center}
\hspace{7,35pt}
\begin{tabular}{|r@{}l|l|r@{\hspace{1ex}}r@{\hspace{0,5ex}}r|c|} \hline
\multicolumn{2}{|c|}{AB}&\multicolumn{1}{c|}{Webs}&\multicolumn{1}{c|}{b}&\multicolumn{1}{c|}{l}&\multicolumn{1}{c|}{n}&s\\
\hline
$9$&$_{9}$&$-2879,-8627$&$3$&$10$&$4$&r\\
$9$&$_{10}$&$16482^{\rm s}$&$4$&$11$&$1$&r\\
$9$&$_{11}$&$-963$&$4$&$11$&$2$&r\\
$9$&$_{12}$&$-295$&$5$&$12$&$2$&r\\
$9$&$_{13}$&$16162,16194$&$4$&$11$&$4$&r\\
$9$&$_{14}$&$127794$&$5$&$14$&$2$&r\\
$9$&$_{15}$&$-279$&$5$&$12$&$2$&r\\
$9$&$_{16}$&$-7763$&$3$&$10$&$2$&r\\
$9$&$_{17}$&$33981^{\rm s},-153066^{\rm s}$&$4$&$13$&$2$&r\\
$9$&$_{18}$&$-859$&$4$&$11$&$2$&r\\
$9$&$_{19}$&$31605,-128766$&$5$&$14$&$4$&r\\
$9$&$_{20}$&$-899,-2687$&$4$&$11$&$4$&r\\
$9$&$_{21}$&$14322$&$5$&$12$&$2$&r\\
$9$&$_{22}$&$33693,-152885$&$4$&$13$&$4$&r\\
$9$&$_{23}$&$-20858^{\rm s},22882^{\rm s}$&$4$&$11$&$2$&r\\
$9$&$_{24}$&$1235,3695$&$4$&$11$&$4$&r\\
$9$&$_{25}$&$-25407,29847,33279,$&&&&\\
&&$-32897,-69055$&$5$&$14$&$10$&r\\
$9$&$_{26}$&$-33677,152890,$&&&&\\
&&$152922$&$4$&$13$&$6$&r\\
$9$&$_{27}$&$-915$&$4$&$11$&$2$&r\\
$9$&$_{28}$&$7576,-7589,-23104$&$4$&$11$&$6$&r\\
$9$&$_{29}$&$-691971$&$4$&$15$&$2$&r\\
$9$&$_{30}$&$-8169$&$4$&$11$&$2$&r\\
$9$&$_{31}$&${\it-2354}^{\rm s}$&$4$&$9$&$1$&r\\
$9$&$_{32}$&$-3771^{\rm n},-108717^{\rm n}$&$4$&$13$&$2$&n\\
$9$&$_{33}$&$23739^{\rm n},-69148^{\rm n},$&&&&\\
&&$101459^{\rm n}$&$4$&$13$&$3$&n\\
$9$&$_{34}$&$-1362201,-1976361,$&&&&\\
&&$12383186$&$4$&$17$&$6$&r\\
$9$&$_{35}$&$436350,3427966$&$5$&$16$&$4$&r\\
$9$&$_{36}$&$-8051$&$4$&$11$&$2$&r\\
$9$&$_{37}$&$94809^{\rm s},-168523,$&&&&\\
&&$218815,-3474018^{\rm s}$&$5$&$16$&$6$&r\\
$9$&$_{38}$&$145570$&$4$&$13$&$2$&r\\
$9$&$_{39}$&$500919,3480870$&$5$&$16$&$4$&r\\
$9$&$_{40}$&$-5915517$&$4$&$17$&$2$&r\\
$9$&$_{41}$&$12114242,-15094041,$&&&&\\
&&$-35669036$&$5$&$18$&$6$&r\\
$9$&$_{42}$&$-3191$&$4$&$11$&$2$&r\\
$9$&$_{43}$&$-3197,3255$&$4$&$11$&$4$&r\\
$9$&$_{44}$&$-8057,8115$&$4$&$11$&$4$&r\\
$9$&$_{45}$&$-3309$&$4$&$11$&$2$&r\\
$9$&$_{46}$&$9627$&$4$&$13$&$2$&r\\
$9$&$_{47}$&$-9573$&$4$&$13$&$2$&r\\
$9$&$_{48}$&$145522$&$4$&$13$&$2$&r\\
$9$&$_{49}$&$-62113,147666$&$4$&$13$&$4$&r\\
\hline
\end{tabular}
%\end{center}
}}
%\end{table}

\skipline
\skipline
\noindent
Minimality has been checked by calculating HOMFLY- and Kauffman-polynomials
for all shorter woven braids; braid indices can be calculated using the MFW-inequality
(see \cite{Jo}).

The column ``AB'' lists the Alexander-Briggs notations, ``b'',``l'' and ``n'' give the braid indices, 
the lengths and the number of minimal webs for the corresponding knot, respectively.
The column 's', taken out of \cite{BZ}, contains information about knot symmetries
in Conway's notation:
\skipsline
\begin{center}
\begin{tabular}{l|c|c}
 & amphicheiral & non-amphicheiral \\
 \hline
 invertible & f & r \\
 \hline
 non-invertible & i & n \\
\end{tabular}
\end{center}
\skipsline
The tuples $(e_1,c_u,\ldots,c_v)$ of the minimal webs are listed in the second column, encoded as integers
$e_1(c_u\!-\!1 + 2c_{u+1} + 6c_{u + 2} + \cdots + 2\cdot 3^{v-u-1}c_v)$.
If the woven braids $\omega$ and $\inv(\omega)$ represent the same oriented knot,
we list only the absolutely smaller of the two corresponding integers.
So most of the integers represent two webs; exceptions are marked by
the following superscripts:
\begin{mylist}{}
\item[n:]
The knot is non-invertible.
\item[s:]
$\omega$ is symmetric, showing that the knot is strongly invertible.
\item[a:]
$\omega$ is antisymmetric, showing the strong negative amphicheirality.
\end{mylist}
\noindent
An entry of the form $\pm j$ indicates that $j$ and $-j$ appear (i.e.~amphicheirality).
Integers in italics represent webs with alternating tight words.
\skipsline
In table \ref{knottable} we get, on average, 2.9 minimal webs per knot type,
having an average tight word length of 11.5.
Symmetric woven braids occur for twelve knots and all seven amphicheiral
knots in the table have antisymmetric minimal woven braids.
The author also found antisymmetric web-representatives for the amphicheiral
ten-crossing knots, being minimal except in the case of $10_{81}$
(a list of minimal webs for all ten-crossing knots is available on request).
This demonstrates that all amphicheiral knots with at most ten crossings are
strongly negative amphicheiral.

It will be the task of a subsequent article to show that any strongly negative
amphicheiral knot can be presented as closure of some antisymmetric woven braid.
\skipsline

\renewcommand{\section}[1]{\noindent\tenbf}

\end{document}